\documentclass{article}

\usepackage[utf8]{inputenc}
\usepackage[T1]{fontenc}
\usepackage{lmodern}
\usepackage{amsmath, amssymb}
\usepackage{graphicx}
\usepackage{hyperref}
\usepackage{geometry}
\usepackage{pgfplots}
\usepgfplotslibrary{groupplots}
\usepackage{subcaption}
\usepackage{siunitx}
\usepackage{placeins}

\usepackage{biblatex}
\usepackage{algorithm}
\usepackage{algpseudocode}

\usepackage{parskip}

\title{Fluid boundary conditions in kinetic-diffusion Monte Carlo}
\author{Thijs Steel \and Vince Maes \and Giovanni Samaey}

\definecolor{ferngreen}{HTML}{56641a}
\definecolor{perfumepurple}{HTML}{c0affb}
\definecolor{apricotorange}{HTML}{e6a176}
\definecolor{orientblue}{HTML}{00678a}
\definecolor{winered}{HTML}{984464}
\definecolor{downygreen}{HTML}{5eccab}

\begin{document}




\maketitle

\begin{abstract}
    The Kinetic-Diffusion Monte Carlo (KDMC) method is a powerful tool for simulating neutral particles in fusion reactors. It is a hybrid fluid-kinetic method that is significantly faster than pure kinetic methods at the cost of a small bias due to fluid approximations. Unfortunately, when simulating particles close to a boundary, it needs to switch to a purely kinetic method, which is significantly slower. In this paper, we will extend the method so that it can accurately take boundary conditions into account without switching to a purely kinetic method. Experiments show that this extension can lead to a speedup of up to 500 times compared to a KDMC method that switches to a purely kinetic method, while not sacrificing too much accuracy.
\end{abstract}

\section{Introduction}

A key challenge in the design of fusion reactors is the simulation of the neutral particles. These are typically modelled using a linear kinetic equation \cite{lapeyre1998methodes}:
\begin{align}\label{eq:boltzmann}
    \underbrace{\partial_t f(x,v,t)}_{\text{transient}} 
    + \underbrace{v \cdot \nabla f(x,v,t)}_{\text{transport}} &= \underbrace{S(x,v,t)}_{\text{sources}} 
    - \underbrace{R_i(x,t) f(x,v,t)}_{\text{ionization}} \notag \\
    &+ \underbrace{R_{cx}(x,t)\left(M(v|x,t)\int f(x,v',t)dv' - f(x,v,t) \right)}_{\text{charge exchange}},
\end{align}
where $f(x,v, t)$ is the density of the neutrals with velocity $v$ at time $t$ with a given position $x$; $S(x, v, t)$ is the source term; $R_i$ is the rate at which particles at a given position $x$ are ionized; $R_{cx}$ is the rate at which particles at a given position $x$ undergo charge exchange; and $M(v|x,t)$ is the post-collisional velocity distribution. We assume that $M(v|x,t)$ is a Gaussian distribution with mean $\nu_p$ and variance $\sigma_p^2$. For the sake of simplicity, we assume that the ionization term is negligible and that there is no explicit time dependence in the coefficients.

The simulation of the neutral particles is closely coupled to the simulation of the charged particles. The source term, the ionization rate, the charge exchange rate, as well as the parameters $\nu_p$ and $\sigma_p$ are all determined by the charged particles. Conversely, the charged particles are also affected by the neutral particles. This leads to iterative schemes where the neutral and charged particles are updated in an alternating fashion until convergence. Schemes like this can be found in, e.g., SOLPS-ITER \cite{wiesen2015new}.

Equation \eqref{eq:boltzmann} is high-dimensional and a common approach to solve it is to use a Monte Carlo method \cite{lapeyre1998methodes}. The idea is to sample a finite number of particles from their initial distribution. Each particle is then traced through the simulation domain, where it moves at a constant velocity until it collides with the background plasma. The collision rate is determined by the ionization and charge exchange rates, which are computed based on the current state of the charged particles. When a particle collides, it either changes its velocity according to the post-collisional distribution $M(v|x,t)$ or is removed from the simulation if it is ionized. The Monte Carlo method is highly parallelizable, but it can be computationally expensive, especially in high-collisional regions where many collisions need to be simulated.

Another common approach is to use a fluid approximation of the kinetic equation \cite{reiter2005eirene}. The idea is to approximate the high-dimensional distribution $f(x,v,t)$ by a combination of low-dimensional functions: density $\rho(x,t)$, mean velocity $\nu(x,t)$, and temperature $T(x,t)$. The kinetic equation can then be rewritten as a set of partial differential equations (PDEs) for these low-dimensional functions. Because of the reduced dimensionality, it is now feasible to solve the PDEs using standard numerical methods, such as finite volumes, which is much cheaper than the Monte Carlo method. The major downside is that the fluid approximation is only valid in high-collisional regions.

In this paper, we focus on a hybrid method that combines the two approaches: kinetic-diffusion Monte Carlo (KDMC) \cite{mortier2022kinetic}. We note that other types of hybrid methods exist, each with their own advantages and disadvantages \cite{borodin2022fluid}. KDMC is a Monte Carlo method where the particles behave more kinetically in low-collisional regions, as they would in a standard Monte Carlo simulation, and more diffusively in high-collisional regions, as they would in a fluid simulation. However, the original formulation of KDMC cannot take boundary conditions into account during the diffusive steps and must switch to a fully kinetic simulation when the particle is close to the boundary. We show that this is not necessary and that it is possible to simulate the diffusive steps accurately close to the boundary. This makes KDMC simulations more efficient and can even make them more accurate.

The remainder of this text is organized as follows. First in Section~\ref{sec:kdmc}, we give a brief overview of the KDMC method. Next, in Section~\ref{sec:boundary_conditions}, we show how to modify the diffusive step in KDMC to make it valid in the presence of boundaries. Finally, in Section~\ref{sec:experiments}, we present some numerical results that demonstrate the effectiveness of the proposed method.

\section{Kinetic-diffusion Monte Carlo}\label{sec:kdmc}

In this section, we give a brief overview of the KDMC method. For a more detailed description, we refer to the original paper \cite{mortier2022kinetic}.

In a standard kinetic Monte Carlo simulation, the particles move with a constant velocity until they collide with the background plasma or the boundary. The time between collisions is exponentially distributed with parameters determined by the charge exchange rate\footnote{In heterogeneous backgrounds, this involves integrating over the charge exchange rate along the trajectory of the particle.}. If the charge exchange rate is high, the particles undergo many collisions. Each of these collisions must be simulated, which can be computationally expensive.

In the fluid limit, we can approximate the kinetic motion with a random walk. In 1D, this means that the particles move according to the following equation:
\begin{equation}\label{eq:diffusion_sde}
    dX_t = \nu(X_t) dt + \sqrt{2D(X_t)} dW_t,
\end{equation}
where $W(t)$ is a standard Wiener process and the fluid coefficients are given by:
\begin{equation}\label{eq:fluid_drift}
    \nu(X_t) = \nu_p(X_t) + \sigma_p^2 \frac{\partial}{\partial x}\frac{1}{R_{cx}},
\end{equation}
for the drift and
\begin{equation}\label{eq:fluid_diffusion}
    D(X_t) = \frac{\sigma_p^2}{R_{cx}},
\end{equation}
for the diffusion. This SDE can then be simulated using a method like Euler-Maruyama. Because it does not need to simulate individual collisions, a single time step of Euler-Maruyama can simulate many collisions at once, which makes it much faster than a standard kinetic Monte Carlo simulation. Unfortunately, the fluid approximation is only valid in high-collisional regions.

KDMC combines kinetic and diffusive motion in a clever way. Each time step begins with a kinetic simulation until the particle collides with the background plasma. This is also called the kinetic step. The rest of the time step is then simulated using a random walk, which is called the diffusive step. The cleverness of this formulation becomes apparent when we consider both low- and high-collisional regions. In low-collisional regions, the time between collisions is long, so the kinetic step takes up most of the time step, which is exactly what we want because the fluid approximation is not valid there. In high-collisional regions, the time between collisions is short, so the kinetic step takes up only a small part of the time step, which is again exactly what we want because the fluid approximation is valid there and much cheaper to simulate. In other words, KDMC automatically switches between kinetic and diffusive motion depending on the collisionality of the region.

During the diffusive step, KDMC does not use the standard fluid coefficients from \eqref{eq:fluid_drift} and \eqref{eq:fluid_diffusion}. Instead, it uses a modified version of the coefficients so that the diffusive step has the same mean and variance as in a kinetic motion. We refer to \cite{mortier2022kinetic}, equations (3.5) and (3.6) for the full expressions of the coefficients. In this text, we will not use these expressions, but we will use an improved formulation of KDMC, called KDKMC \cite{mortier2020kinetic}. In KDKMC, we do another small free flight step after the diffusive step. This improves the accuracy of the method and (most importantly for our purposes), it simplifies the equations for the coefficients. The drift and diffusion coefficients are then given by:
\begin{equation}
    \nu(X_t) = \nu_p + \sigma_p^2\frac{\partial}{\partial x}\frac{1}{R_{cx}},
\end{equation}
\begin{equation}
    D(X_t) = \frac{\sigma_p^2}{R_{cx}^2 \theta}(2 e^{-R_{cx}\theta} + R_{cx}\theta + R_{cx}\theta e^{-R_{cx}\theta} - 2),
\end{equation}
where $\theta$ is the time step of the diffusive step.

Finally, we also mention other extensions of KDMC. KDMC has been extended with fluid estimators \cite{mortier2022estimation}, which is essential for estimating quantities of interest. KDMC has also been incorporated into a multi-level framework \cite{mortier2022multilevel}, which can significantly speed up the simulation.

\section{Boundary conditions in the diffusive step}\label{sec:boundary_conditions}

During the diffusive step, the kinetic motion of the particle is approximated by a random walk. In a homogeneous background without boundaries, the drift and diffusion coefficient are constant and the SDE can be solved exactly. In this case, the solution is given by:
\begin{equation}
    X_t = X_0 + \nu t + \sqrt{2D t} W(t).
\end{equation}
This also matches the result of the Euler-Maruyama method. However, this exact solution is not valid in the presence of boundaries (or in heterogeneous backgrounds). In this section, we derive a way to simulate the SDE in the presence of boundaries.

First, in Subsection~\ref{subsec:solution_sde_boundaries}, we derive an analytical solution of the probability density function of the particle's position that is valid in the presence of boundaries for a 1D equation with homogeneous plasma background. Second, in Subsection~\ref{subsec:sampling_pdf}, we show how we can efficiently sample from this probability density function.

\subsection{Solution of the SDE in the presence of boundaries}\label{subsec:solution_sde_boundaries}

The goal of the subsection is to derive an analytical solution of the SDE \eqref{eq:diffusion_sde} with homogeneous parameters that takes boundaries into account. To find such an analytical solution, we switch from the SDE form to the Fokker-Planck equation associated with the SDE. The Fokker-Planck equation is given by:
\begin{equation}
    \frac{\partial p}{\partial t}(x,t) = -\nu \frac{\partial p}{\partial x}((x,t)) + D \frac{\partial^2 p}{\partial x^2}((x,t)),
\end{equation}
subject to the Robin boundary condition:
\begin{equation}
    \alpha p(L,t) + \beta \frac{\partial p}{\partial x}(L,t) = 0,
\end{equation}
and initial condition:
\begin{equation}
    p(x, 0) = \delta(x - x_0).
\end{equation}
The simulation domain is given by $]-\infty, L]$. A similar PDE is solved in \cite{sommerfeld1949partial} and we follow their approach, modified where necessary, to find the following solution:
\begin{equation}\label{eq:solution_p}
    p(x, t) = \underbrace{U(x,t,x_0)}_{p_1(x,t)} + \underbrace{e^{\frac{\nu}{D}(L - x_0)} U(x,t,x_R)}_{p_2(x,t)} + \underbrace{2(\frac{\alpha}{\beta} + \frac{\nu}{2D}) e^{\frac{\nu}{D}(L - x_0)} \int_{x_R}^{\infty} e^{\frac{\alpha}{\beta}(\eta - x_R)} U(x,t,\eta) d\eta}_{p_3(x,t)},
\end{equation}
where $x_R = 2L - x_0$ and $U(x, t; \eta)$ is the free solution centered around $\eta$:
\begin{equation}
    U(x, t, \eta) = \frac{1}{\sqrt{4 \pi D t}} e^{-\frac{(x - \eta - \nu t)^2}{4Dt}}.
\end{equation}
$p_1(x,t)$ 
This solution is only valid if $\beta$ is nonzero. If $\beta = 0$, we have a purely absorbing boundary and the solution simplifies to:
\begin{equation}\label{eq:solution_p_absorbing}
    p(x, t) = U(x,t,x_0) - e^{\frac{\nu}{D}(L - x_0)} U(x,t,x_R).
\end{equation}

\subsection{Sampling from the probability density function}\label{subsec:sampling_pdf}

Now that we have found an analytical solution for the probability density function (pdf) $p(x, t)$ of the particle's position in the presence of boundaries, we need to sample from this probability density function to use it in the diffusive step of KDMC. We dicuss two sampling methods, a basic and a more efficient version.

\subsubsection{Basic sampling method}

First, we focus on $p_1(x,t)$ and $p_2(x,t)$. Together, these terms form a mixture of two Gaussians and we can sample from them by first drawing a random number to decide which Gaussian to sample from and then sampling from the selected Gaussian. We note that we only consider samples that lie within the domain. Many methods exist to sample from a truncated normal; we will use a method described C. P. Robert \cite{robert1995simulation}. This is an accept-reject method with exponential proposals that is efficient enough for our purposes.

Next, we focus on $p_3(x,t)$, which is negative for practical values of $\alpha$ and $\beta$. We use an accept-reject scheme to take it into account. We propose samples from $p_1(x,t) + p_2(x,t)$ and accept them with the following rate: $\frac{p(x,t)}{p_1(x,t) + p_2(x,t)}$. This is allowed because $p_1(x,t) + p_2(x,t) > p(x,t)$ for all $x$ and $t$ (assuming $p_3(x,t)$ is negative).

Finally, we must take into account that $p(x,t)$ is not (necessarily) normalized. Because of absorbtion at the wall, it is possible that $Q = \int_{-\infty}^L p(x,t) dx < 1$. To take this into account, we just need to multiply the weight of all the particles by $Q$.

\subsubsection{More efficient sampling method}

The basic sampling method has one major disadvantage: it always takes the boundary into account, even for particles that do not cross the boundary. This is not really a problem in 1D, but it could be a considerable bottleneck in 2D/3D, where there could be a large amount of boundaries to consider.

We propose the following alternative method. First, we sample from $p_1(x,t)$ without considering the boundary. If it is within the boundary, we accept is as a valid sample and do not adjust its weight. If it is outside the boundary, we sample from $p_2(x,t) + p_3(x,t)$ using proposals from $p_2(x,t)$ and accepting them with the following rate: $\frac{p_2(x,t) + p_3(x,t)}{p_2(x,t)}$. We then multiply the weight of these particles by $\frac{Q_1}{Q_2}$, where $Q_1 = \int_{-\infty}^L p_2(x,t) + p_3(x,t) dx$ and $Q_2 = \int_{L}^\infty p_1(x,t) dx$. We note that this scheme only works if $p_2(x,t) + p_3(x,t) > 0$. To make it work in all cases, we would have to introduce negative weights.

\section{Experiments}\label{sec:experiments}

In this section, we present some numerical experiments to demonstrate the effectiveness of the proposed method. We have implemented the different methods in a newly developed C++ package called NEPTUNE. It is publically available at \url{https://gitlab.kuleuven.be/numa/software/neptune-mc}. We compare the following methods:
\begin{itemize}
    \item Kinetic reference: a standard kinetic Monte Carlo simulations
    \item Fluid model: a simple advection-diffusion fluid model using the coefficients described in \eqref{eq:fluid_drift} and \eqref{eq:fluid_diffusion}.
    \item KDMC\_Kin: The base algorithm is KDKMC as described in \cite{mortier2020kinetic}. When a particle crosses the boundary during the diffusive step, or when the mean of the diffusive step is less than 2 standard deviations away from the boundary, the particle is simulated kinetically for the rest of the time step.
    \item KDMC\_Fluid: The base algorithm is KDKMC as described in \cite{mortier2020kinetic}. It never switches to a fully kinetic simulation, but instead uses the modified diffusive step described in Section~\ref{sec:boundary_conditions}.
\end{itemize}

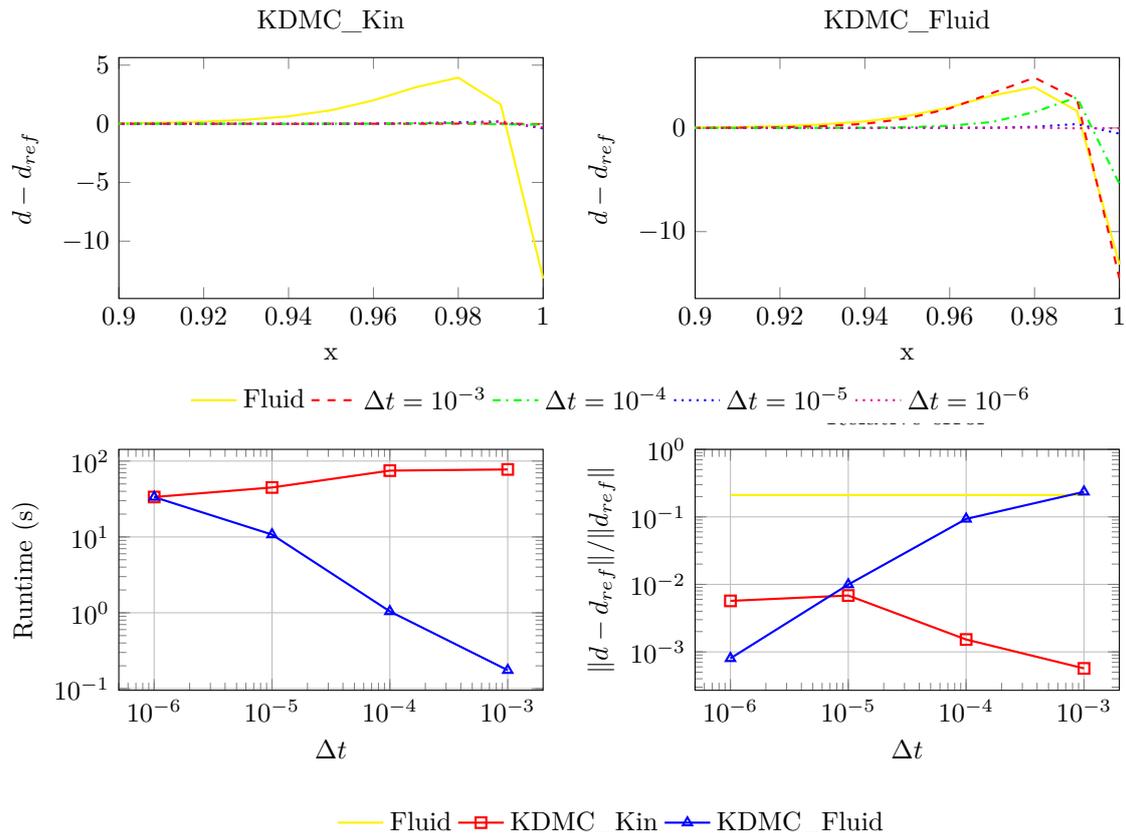
\begin{figure}[htb!]
    \centering
    \begin{tikzpicture}
        
\begin{groupplot}[group style={
        {group size=2 by 2, horizontal sep=2cm, vertical sep=2cm}},name=t1, height=0.3\textwidth, width=0.45\textwidth,legend style={
        transpose legend,
        legend columns=0,
        draw=none }]

    \nextgroupplot[title=KDMC\_Kin, legend to name=grouplegend, xmin=0.9, xmax=1, xlabel={x}, ylabel={$d - d_{ref}$}]

            \addplot[yellow, thick, mark=none] table [x=x, y expr={\thisrow{fluid} - \thisrow{ref}}, col sep=comma] {data/1D_reflecting_drift.csv};
            \addlegendentry{Fluid}

            \addplot[red, thick, mark=none, dashed] table [x=x, y expr={\thisrow{kd_old_0.001} - \thisrow{ref}}, col sep=comma] {data/1D_reflecting_drift.csv};
            \addlegendentry{$\Delta t = 10^{-3}$}

            \addplot[green, thick, mark=none, dashdotted] table [x=x, y expr={\thisrow{kd_old_0.0001} - \thisrow{ref}}, col sep=comma] {data/1D_reflecting_drift.csv};
            \addlegendentry{$\Delta t = 10^{-4}$}

            \addplot[blue, thick, mark=none, dotted] table [x=x, y expr={\thisrow{kd_old_1e-05} - \thisrow{ref}}, col sep=comma] {data/1D_reflecting_drift.csv};
            \addlegendentry{$\Delta t = 10^{-5}$}

            \addplot[magenta, thick, mark=none, dotted] table [x=x, y expr={\thisrow{kd_old_1e-06} - \thisrow{ref}}, col sep=comma] {data/1D_reflecting_drift.csv};
            \addlegendentry{$\Delta t = 10^{-6}$}

    \nextgroupplot[title=KDMC\_Fluid, xmin=0.9, xmax=1, xlabel={x}, ylabel={$d - d_{ref}$}]

            \addplot[yellow, thick, mark=none] table [x=x, y expr={\thisrow{fluid} - \thisrow{ref}}, col sep=comma] {data/1D_reflecting_drift.csv};

            \addplot[red, thick, mark=none, dashed] table [x=x, y expr={\thisrow{kd_new_0.001} - \thisrow{ref}}, col sep=comma] {data/1D_reflecting_drift.csv};

            \addplot[green, thick, mark=none, dashdotted] table [x=x, y expr={\thisrow{kd_new_0.0001} - \thisrow{ref}}, col sep=comma] {data/1D_reflecting_drift.csv};

            \addplot[blue, thick, mark=none, dotted] table [x=x, y expr={\thisrow{kd_new_1e-05} - \thisrow{ref}}, col sep=comma] {data/1D_reflecting_drift.csv};

            \addplot[magenta, thick, mark=none, dotted] table [x=x, y expr={\thisrow{kd_new_1e-06} - \thisrow{ref}}, col sep=comma] {data/1D_reflecting_drift.csv};




    \nextgroupplot[ymode=log, xmode=log, xlabel={$\Delta t$}, ylabel={Runtime (s)}, grid=major]

            \addplot[thick, red, mark=square] table [x=dt, y expr={\thisrow{runtime_old}}, col sep=comma] {data/1D_reflecting_drift_runtime.csv};

            \addplot[thick, blue, mark=triangle] table [x=dt, y expr={\thisrow{runtime_new}}, col sep=comma] {data/1D_reflecting_drift_runtime.csv};

    \nextgroupplot[title=Relative error, ymode=log, xmode=log,xlabel={$\Delta t$}, ylabel={$\|d - d_{ref}\|/\|d_{ref}\|$}, log basis y=10, legend to name=grouplegend2, grid=major, ymax=1]

            \addplot[thick, yellow, mark=none] table [x=dt, y expr={\thisrow{error_fluid}}, col sep=comma] {data/1D_reflecting_drift_runtime.csv};
            \addlegendentry{Fluid}

            \addplot[thick, red, mark=square] table [x=dt, y expr={\thisrow{error_old}}, col sep=comma] {data/1D_reflecting_drift_runtime.csv};
            \addlegendentry{KDMC\_Kin}

            \addplot[thick, blue, mark=triangle] table [x=dt, y expr={\thisrow{error_new}}, col sep=comma] {data/1D_reflecting_drift_runtime.csv};
            \addlegendentry{KDMC\_Fluid}
\end{groupplot}

    \node at (6.5,0)
    [below, yshift=-2\pgfkeysvalueof{/pgfplots/every axis title shift}]
    {\ref{grouplegend}};

    \node at (6.5,0)
    [below, yshift=-18\pgfkeysvalueof{/pgfplots/every axis title shift}]
    {\ref{grouplegend2}};
    

\end{tikzpicture}

\caption{Comparison of KDMC\_Kin and KDMC\_Fluid for a 1D reflecting boundary with drift. In the top left, the difference with the kinetic reference is shown for the density estimated by KDMC\_Kin and a fully fluid simulation. The top right plot shows the same for KDMC\_Fluid. The bottom left plot shows the speedup of both KDMC methods compared to the fully kinetic simulation for different time steps. The bottom right plot shows the relative normwise error of the fluid and both KDMC methods compared to the fully kinetic simulation for different time steps. }

\label{fig:1D_reflecting_boundary_drift}
\end{figure}

We will simulate a 1D domain with either reflecting boundaries. The domain is given by $[0, 1\unit{\metre}]$, discretized into 101 cells and the initial position of the particle is $x_0 = 0.98\unit{\metre}$ with an initial velocity drawn from the plasma background, $\nu_p = 100\unit{\metre \per \second}$, $R_{cx} = 1.0e7\unit{\per \second}$, and $\sigma_p^2 = 1.0e7\unit{\metre^2 \per \second^2}$. We simulate a total of $10^6$ particles until a final time of $t_f = 0.01\unit{\second}$. The KDMC time step $\Delta t$ is varied between $10^{-6}\unit{\second}$ and $10^{-3}\unit{\second}$. The results are shown in Figure \ref{fig:1D_reflecting_boundary_drift}.

The following things are noteworthy about the results:
\begin{itemize}
    \item It should be immediately clear that KDMC is significantly faster when using fluid boundary conditions. The speedup is especially noticeable for larger time steps, where KDMC\_Fluid is up to 500 times as fast as KDMC\_Kin. This speedup can be explained by the fact that KDMC\_Fluid does not need to switch to a fully kinetic simulation when the particle is close to the boundary. For large time steps, up to 100\% of the particle trajectories are simulated kinetically in KDMC\_Kin, which explains the small speedup KDMC\_Kin achieves over the kinetic reference.
    \item KDMC\_Fluid is generally accurate, but despite being based on an exact analytical solution, it is not a perfect match with the kinetic reference solution. This is because the fluid model is only an approximation of the kinetic equations. This is especially visible in for $\Delta t = 10^{-3}$ where the error is very similar to that of the fluid model.
    \item KDMC\_Kin shows some strange behavior for increasing timesteps. We would expect both the error and the speedup to increase with larger timesteps, but the opposite is true. This is because larger timesteps lead to more particles crossing the boundary during a diffusive step, and thus more particles being simulated kinetically. As kinetic simulations are more accurate and expensive, this explains the behavior. This effect can be severe; for $\Delta t = 10^{-3}$, 100\% of the particle trajectories are simulated kinetically, while for $\Delta t = 10^{-6}$, only 19\% of the particle trajectories are simulated kinetically.
\end{itemize}

In KDMC, we expect the solution to be close to that of the fluid model for large time steps and close to the kinetic reference solution for small time steps. KDMC\_Fluid exhibits this behavior, but KDMC\_Kin does not. In fact, as the time step increases, KDMC\_Kin is a better match with the kinetic reference solution than it was for smaller time steps. This is because for these large time steps, the chance that a particle will cross the boundary and thus need to be simulated kinetically is very high.

\section{Conclusion}
In this paper, we have shown how to modify the diffusive step in KDMC to make it valid in the presence of boundaries. We derived an analytical solution for the probability density function of the particle's position in the presence of boundaries and showed how to sample from it efficiently. Experiments show that the new method is significantly faster than the original method.

In future work, we will extend these boundary techniques to higher dimensions and we will improve the accuracy of the SDE simulation in heterogeneous backgrounds.

\section*{Author contributions}

\textbf{Thijs Steel}: Conceptualization (lead), Writing software, Writing of original draft. \textbf{Vince Maes}: Conceptualization(support), Review and editing. \textbf{Giovanni Samaey}: Conceptualization(support), Review and editing, Supervision, Funding acquisition.

\section*{Acknowledgments}

The authors are grateful to Wouter Dekeyser for many useful discussions.

\section*{Financial disclosure}

This work has been carried out within the framework of the EUROfusion
Consortium, funded by the European Union via the Euratom Research and
Training Programme (Grant Agreement No 101052200 — EUROfusion). Views and
opinions expressed are, however, those of the author(s) only and do not necessarily
reflect those of the European Union or the European Commission. Neither the
European Union nor the European Commission can be held responsible for them.

Part of this research was funded by the Research Foundation Flanders (FWO) under grant 3E211335. V. Maes is a SB PhD fellow of the Research Foundation
Flanders (FWO), funded by grant 1S64723N.

\section*{Conflict of interest}

The authors declare that there is no conflict of interests.

\printbibliography


\end{document}